\documentclass[reqno,12pt]{amsart}
\usepackage{amsmath,amssymb,tikz,stix2}
\usepackage[colorlinks=true,
linkcolor=red,
citecolor=green,
urlcolor=green]{hyperref}

\newtheorem{thm}{Theorem}
\newtheorem{lem}[thm]{Lemma}
\newtheorem{cor}[thm]{Corollary}
\newtheorem{prop}[thm]{Proposition}

\theoremstyle{definition}

\theoremstyle{remark}

\DeclareMathOperator{\essinf}{ess\,inf}
\DeclareMathOperator{\esssup}{ess\,sup}

\begin{document}
\emergencystretch=3em

\title[connectedness of the spectra of self-adjoint
Toeplitz operators]
{On the connectedness of the spectra of self-adjoint
Toeplitz operators on the sphere}

 \author{Shibo Xu}
\address[Shibo Xu]{School of Mathematical Sciences, Dalian University of Technology, 
Dalian, Liaoning 116024, P. R. China}
\email{shiboxu98@163.com}

\author{Congquan Yan}
\address[Congquan Yan]{College of Data Science, Jiaxing University, Jiaxing, Zhejiang 314001, P. R. China}
\email{yancongquan@zjxu.edu.cn}

\author{Tao Yu}
\address[Tao Yu]{School of Mathematical Sciences, Dalian University of Technology,
Dalian, Liaoning 116024, P. R. China}
\email{tyu@dlut.edu.cn}

\subjclass[2020]{Primary 47A10, 47B35.}
\keywords{Toeplitz operator, Hardy space, spectrum, Fredholm index.}

\begin{abstract}
Let $H^2(\partial\mathbb B_n)$ be the Hardy space of the unit sphere
in $\mathbb C^n$, $n\geq 2$. We prove that every self-adjoint
Toeplitz operator on $H^2(\partial\mathbb B_n)$ has connected
spectrum. More precisely, for every real-valued
$f\in L^\infty(\partial\mathbb B_n)$, both the spectrum and
the Fredholm essential spectrum of $T_f$ coincide with
the interval $[\essinf f,\esssup f]$.
\end{abstract}

\maketitle

\section{Introduction}
The connectedness of the spectra of Toeplitz operators is a
classical problem in operator theory. On the Hardy space
$H^2(\mathbb T)$ of the unit circle, Hartman and Wintner proved that the spectrum of a self-adjoint Toeplitz operator $T_f$ on the Hardy space of the unit circle
coincides with the interval
$[\operatorname*{ess\,inf} f,\operatorname*{ess\,sup} f]$. In 1964, Widom \cite{Widom}, answering a question of Halmos,
proved that the spectrum of every Toeplitz operator with a bounded
measurable symbol on $\mathbb T$ is connected. Subsequently, Douglas
established the connectedness of the essential spectrum of such
operators; see \cite{Douglas}. The extent to which these spectral
properties persist in several complex variables has been a natural
question in the theory of Toeplitz operators on function spaces.

Let $\mathbb B_n$ be the open unit ball in $\mathbb C^n$ 
and let $\partial\mathbb B_n$ be the unit sphere equipped 
with the normalized surface measure $d\sigma$. The Hardy 
space $H^2(\partial\mathbb B_n)$ is the closed subspace of 
$L^2(\partial\mathbb B_n,d\sigma)$ generated by the 
holomorphic polynomials. Let $P$ denote the orthogonal 
projection of $L^2(\partial\mathbb B_n,d\sigma)$ onto 
$H^2(\partial\mathbb B_n)$. Given a symbol 
$f\in L^\infty(\partial\mathbb B_n)$, the Toeplitz operator 
$T_f$ on $H^2(\partial\mathbb B_n)$ is defined by
\[
T_f=PM_f\big|_{H^2(\partial\mathbb B_n)},
\]
where $M_f$ is the multiplication operator induced by $f$ 
on $L^2(\partial\mathbb B_n,d\sigma)$.
We denote by $\sigma(T)$ and $\sigma_{\mathrm e}(T)$
the spectrum and the Fredholm essential spectrum of
a bounded operator $T$, respectively. Thus
\[
\sigma_{\mathrm e}(T)
=\{\lambda\in\mathbb C:
T-\lambda I\text{ is not Fredholm}\}.
\]

In their seminal work on Toeplitz operators in several complex
variables, Davie and Jewell \cite{DavieJewell} established the
spectral inclusion
\[
\operatorname{ess\,ran} f\subseteq \sigma(T_f),
\qquad f\in L^\infty(\partial\mathbb B_n).
\]
They also showed that, in contrast to the one variable
case, neither the spectrum nor the essential spectrum of a
Toeplitz operator on $H^2(\partial\mathbb B_n)$, $n\geq2$, is
necessarily connected. The counterexamples, however, do not
involve self-adjoint Toeplitz operators. This led Davie and Jewell
\cite[p.~359]{DavieJewell} in 1977 to conjecture that
\[
\sigma(T_f)=[\essinf f,\esssup f]
\]
for every real-valued $f\in L^\infty(\partial\mathbb B_n)$.

An affirmative answer to the conjecture for a special class of
symbols was recently obtained by Xia and Yan \cite{XiaYan}.
More precisely, they proved that the conjecture holds
whenever $f\in L^\infty(\partial\mathbb B_n)$ is real-valued and
invariant under coordinatewise rotations, that is,
\[
f(\tau_1z_1,\ldots,\tau_nz_n)
=f(z_1,\ldots,z_n)
\]
almost everywhere on $\partial\mathbb B_n$, for every
$(\tau_1,\ldots,\tau_n)\in\mathbb T^n$.

The rotational invariance of the symbol ensures that
$T_f$ is diagonal with respect to the standard
orthonormal basis
$\{e_\alpha:\alpha\in\mathbb N_0^n\}$
of $H^2(\partial\mathbb B_n)$,
where $\mathbb N_0=\{0,1,2,\ldots\}$. Consequently,
\[
\sigma(T_f)
=\overline{\left\{
\langle T_fe_\alpha,e_\alpha\rangle:
\alpha\in\mathbb N_0^n
\right\}}.
\]
The proof of Xia and Yan relies on two ingredients:
the approximation of the essential infimum and essential supremum
of $f$ by eigenvalues indexed by multi-indices of sufficiently large
total degree, via multivariate Bernstein polynomials,
and estimates for the differences between eigenvalues
corresponding to neighboring multi-indices.
By constructing suitable chains of such eigenvalues,
they establish their density in $[\essinf f,\esssup f]$.

This approach depends essentially on the diagonal structure
of $T_f$ and does not directly extend to arbitrary real-valued
symbols.

The purpose of the present paper is to establish the
Davie--Jewell conjecture in its full generality. Our main result
also identifies the essential spectrum of every self-adjoint
Toeplitz operator on the spherical Hardy space.

\begin{thm}\label{thm}
Let $n\geq 2$ and let
$f\in L^\infty(\partial\mathbb B_n)$ be real-valued.
Then
\begin{equation}\label{eq:main}
\sigma(T_f)=\sigma_{\mathrm e}(T_f)
=\bigl[\operatorname*{ess\,inf}f,
       \operatorname*{ess\,sup}f\bigr].
\end{equation}
\end{thm}

Our proof takes a different approach, combining boundary
limits via automorphisms of the unit ball with
Fredholm index theory. The key observation is that a
hypothetical spectral gap of $T_f$ gives rise to a spectral
projection whose associated Fredholm index, defined through
the Dirac operator of the coordinate multiplication tuple
\cite{Arveson}, is invariant under ball automorphisms.
Boundary limits, however, force this index to take
two distinct values, yielding a contradiction. This argument
requires neither rotational invariance nor additional
regularity of the symbol.

The paper is organized as follows. In Section~\ref{sec:localization}, we establish the compact commutator properties of
the coordinate multiplication operators and develop the boundary
limits argument. Section~\ref{sec:fredholm} is devoted to the Fredholm index
of the associated Dirac operator and its compressions, as well
as the invariance of the index under ball automorphisms.
In Section~\ref{sec:mainproof}, we prove Theorem~\ref{thm}. Finally, in Section~\ref{sec:consequences}, we derive
spectral consequences for Hankel operators with
characteristic function symbols and extend the main result to
weighted Hardy spaces with norms equivalent to the usual Hardy
norm.

\section{Compact commutators and boundary limits}\label{sec:localization}
\subsection{Shift operators and compact commutators}

We first establish some compactness properties of the coordinate
multiplication operators on $H^2(\partial\mathbb B_n)$. For
$1\leq j\leq n$, let
$S_j=M_{z_j}|_{H^2(\partial\mathbb B_n)}$.
Given $h\in L^\infty(\partial\mathbb B_n)$, we denote by
\[
H_h=(I-P)M_h\big|_{H^2(\partial\mathbb B_n)}
:H^2(\partial\mathbb B_n)\longrightarrow
H^2(\partial\mathbb B_n)^\perp
\]
the Hankel operator with symbol $h$.
 For bounded operators $A$ and $B$ on a Hilbert space,
we denote their commutator by $[A,B]=AB-BA$. For Hilbert spaces $\mathcal H_1$ and $\mathcal H_2$,
let $\mathcal B(\mathcal H_1,\mathcal H_2)$ and
$\mathcal K(\mathcal H_1,\mathcal H_2)$ denote the spaces
of bounded and compact operators from $\mathcal H_1$
to $\mathcal H_2$, respectively.

We use the standard multi-index notation
$|\alpha|=\alpha_1+\cdots+\alpha_n$,
$\alpha!=\alpha_1!\cdots\alpha_n!$, and
$z^\alpha=z_1^{\alpha_1}\cdots z_n^{\alpha_n}$
for $\alpha\in\mathbb N_0^n$. Recall that the functions
\begin{equation}\label{eq:basis}
e_\alpha(z)=
\left(\frac{(|\alpha|+n-1)!}{(n-1)!\,\alpha!}\right)^{1/2}
z^\alpha,
\qquad \alpha\in\mathbb N_0^n,
\end{equation}
form an orthonormal basis for $H^2(\partial\mathbb B_n)$.
For each $m\geq 0$, let
\[
H_m=\operatorname{span}\{e_\alpha:|\alpha|=m\}
\]
be the subspace of homogeneous polynomials of degree $m$.
Then we have the orthogonal decomposition
\[
H^2(\partial\mathbb B_n)=\bigoplus_{m=0}^{\infty}H_m.
\]

With respect to the orthonormal basis in \eqref{eq:basis},
the coordinate multiplication operators act as weighted
shifts:
\begin{equation}\label{eq:shiftformula}
S_je_\alpha=
\left(\frac{\alpha_j+1}{|\alpha|+n}\right)^{1/2}
e_{\alpha+\varepsilon_j},
\end{equation}
where $\varepsilon_j$ denotes the multi-index whose $j$th
entry is one and all others are zero.
By the classical compactness result for Toeplitz
semicommutators with continuous symbols
(see \cite{Coburn} and \cite[pp.~22--23]{Zheng1997}), we have
\begin{equation}\label{eq:essential-shift}
\sum_{j=1}^n S_jS_j^*-I
\in\mathcal K(H^2(\partial\mathbb B_n)),
\qquad
[S_j^*,S_i]
\in\mathcal K(H^2(\partial\mathbb B_n))
\end{equation}
for all $1\leq i,j\leq n$.
In particular, the commuting tuple $(S_1,\ldots,S_n)$
is essentially normal.

\begin{lem}\label{lem:hankelcompact}
For $h\in L^\infty(\partial\mathbb B_n)$ and $1\leq j\leq n$,
\begin{equation}\label{eq:ThS}
[T_h,S_j]=H_{\overline z_j}^{\,*}H_h
\in\mathcal K(H^2(\partial\mathbb B_n)).
\end{equation}
\end{lem}

\begin{proof}
Since $\overline z_j\in C(\partial\mathbb B_n)$,
$H_{\overline z_j}$ is compact by \cite[p.~23]{Zheng1997}.
By the standard Toeplitz--Hankel identity and the fact that
$H_{z_j}=0$, we have
\[
[T_h,S_j]
=T_{hz_j}-T_{z_j}T_h
=H_{\overline z_j}^{\,*}H_h.
\]
The assertion follows.
\end{proof}

We next record some standard consequences of the continuous functional calculus that will be used in the sequel.

\begin{lem}\label{lem:funccalc}
Let $\{A_k\}$ be a sequence of self-adjoint operators on a Hilbert
space $\mathcal H$ with $\sup_k\|A_k\|\leq L$. Let
$B\in\mathcal B(\mathcal H)$ and $h\in C([-L,L])$.
Then the following assertions hold.
\begin{enumerate}
\item[(i)] If $A_k\to cI$ strongly, then
$h(A_k)\to h(c)I$ strongly.

\item[(ii)] If $\|[A_k,B]\|\to0$, then
$\|[h(A_k),B]\|\to0$.

\item[(iii)] If $A=A^*$, $\|A\|\leq L$, and $[A,B]$ is compact,
then $[h(A),B]$ is compact.
\end{enumerate}
\end{lem}

\begin{proof}
(i) By \cite[Chapter~IX, Section~5, Exercise~4]{Conway},
$p(A_k)\to p(c)I$ strongly for every polynomial $p$.
The assertion then follows from the continuous functional
calculus and uniform polynomial approximation on $[-L,L]$.

(ii) For every integer $r\geq1$, the identity
\[
[A_k^r,B]
=\sum_{\ell=0}^{r-1}A_k^\ell[A_k,B]A_k^{r-1-\ell}
\]
and the uniform boundedness of $\{A_k\}$ imply that
$\|[p(A_k),B]\|\to0$ for every polynomial $p$.
The assertion follows by uniform polynomial approximation
on $[-L,L]$.

(iii) Since $\mathcal K(\mathcal H)$ is a norm-closed
two-sided ideal in $\mathcal B(\mathcal H)$
(see \cite{Douglas}), the same commutator identity
and uniform polynomial approximation yield
$[h(A),B]\in\mathcal K(\mathcal H)$.
\end{proof}

\subsection{Ball automorphisms and boundary limits}

For $a\in\mathbb B_n$, let $P_a$ be the orthogonal projection
of $\mathbb C^n$ onto $\mathbb Ca$, and set
$P_a^\perp=I-P_a$ and $s_a=(1-|a|^2)^{1/2}$.
Following \cite[Section~2.2]{Rudin}, we define
\begin{equation}\label{eq:autformula}
\varphi_a(z)=
\frac{a-P_az-s_aP_a^\perp z}{1-\langle z,a\rangle},
\qquad
j_a(z)=\frac{(1-|a|^2)^{n/2}}
{(1-\langle z,a\rangle)^n}.
\end{equation}
Here $P_0=0$, so that $\varphi_0(z)=-z$.
Then $\varphi_a$ is an involutive automorphism of
$\mathbb B_n$ interchanging $0$ and $a$.
By the standard formulas for automorphisms of the unit ball, we have
\begin{equation}\label{eq:autidentity}
j_a(j_a\circ\varphi_a)=1,\qquad
d\sigma(\varphi_a(\zeta))
=|j_a(\zeta)|^2\,d\sigma(\zeta).
\end{equation}

\begin{lem}\label{lem:covariance}
For $a\in\mathbb B_n$, define
\[
(U_au)(\zeta)=j_a(\zeta)u(\varphi_a(\zeta)),
\qquad u\in L^2(\partial\mathbb B_n).
\]
Then $U_a$ is a self-adjoint unitary on
$L^2(\partial\mathbb B_n)$ satisfying $U_aP=PU_a$.
Moreover, for every $h\in L^\infty(\partial\mathbb B_n)$,
\begin{equation}\label{eq:covariance}
U_aT_hU_a^*=T_{h\circ\varphi_a},
\qquad
U_aS_jU_a^*
=M_{\varphi_{a,j}}\big|_{H^2(\partial\mathbb B_n)},
\end{equation}
where $\varphi_{a,j}$ denotes the $j$th coordinate
function of $\varphi_a$.
For each $1\le j\le n$, the map
\[
a\longmapsto M_{\varphi_{a,j}}\big|_{H^2(\partial\mathbb B_n)}
\]
is norm-continuous on $\mathbb B_n$.
\end{lem}

\begin{proof}
By \cite[Theorems~3.3.8 and ~2.2.2]{Rudin},
$U_a$ is a unitary operator on $L^2(\partial\mathbb B_n)$
satisfying $U_aP=PU_a$.
Moreover, by \eqref{eq:autidentity}, we have $U_a^2=I$,
and hence $U_a^*=U_a$.
The identities in \eqref{eq:covariance} follow
by compressing
\[
U_aM_hU_a^*=M_{h\circ\varphi_a}
\]
to $H^2(\partial\mathbb B_n)$.

Finally, it follows directly from \eqref{eq:autformula} that
\[
\|\varphi_{a,j}-\varphi_{b,j}\|_\infty\longrightarrow 0
\qquad\text{as }a\to b.
\]
The asserted norm continuity then follows from
\[
\|M_{\varphi_{a,j}}-M_{\varphi_{b,j}}\|
\leq \|\varphi_{a,j}-\varphi_{b,j}\|_\infty.
\]
\end{proof}

The following lemma shows that bounded measurable symbols
admit constant limits under automorphic translations
at almost every boundary point.

\begin{lem}\label{lem:scalar}
Let $g\in L^\infty(\partial\mathbb B_n)$ be real-valued.
Then, for almost every $\zeta\in\partial\mathbb B_n$,
\[
\lim_{r\to1^-}
\|g\circ\varphi_{r\zeta}-g(\zeta)\|_{L^2(\partial\mathbb B_n)}
=0.
\]
\end{lem}

\begin{proof}
We follow the argument of \cite[Lemma~4.11]{Xia2010}.
For $u\in L^1(\partial\mathbb B_n)$, let
$\mathcal P[u]$ denote its invariant Poisson integral, given by
\[
\mathcal P[u](a)
=\int_{\partial\mathbb B_n}
u(\eta)|j_a(\eta)|^2\,d\sigma(\eta),
\qquad a\in\mathbb B_n,
\]
where $j_a$ is the normalized reproducing kernel of
$H^2(\partial\mathbb B_n)$ at $a$.

Since $g,g^2\in L^1(\partial\mathbb B_n)$, it follows from
\cite[Theorem~5.3.1]{Rudin} that almost every
$\zeta\in\partial\mathbb B_n$ is a Lebesgue point of both
$g$ and $g^2$. For each such $\zeta$, by
\cite[Theorem~5.4.8]{Rudin}, we have
\[
\lim_{r\to1^-}\mathcal P[g](r\zeta)=g(\zeta),
\qquad
\lim_{r\to1^-}\mathcal P[g^2](r\zeta)=g(\zeta)^2.
\]
Consequently, the change of variables formula
\eqref{eq:autidentity} yields
\[
\begin{aligned}
\|g\circ\varphi_{r\zeta}-g(\zeta)\|_{L^2}^2
&=\|(g-g(\zeta))j_{r\zeta}\|_{L^2}^2\\
&=\mathcal P[g^2](r\zeta)
-2g(\zeta)\mathcal P[g](r\zeta)
+g(\zeta)^2\\
&\longrightarrow 0
\end{aligned}
\]
as $r\to1^-$.
\end{proof}

We next pass from convergence of symbols to convergence
of the associated operators.

\begin{lem}\label{lem:localcomm}
Let $\{g_k\}$ be a uniformly bounded sequence of real-valued
functions in $L^\infty(\partial\mathbb B_n)$ such that
$g_k\to c$ in $L^2(\partial\mathbb B_n)$ for some $c\in\mathbb R$.
Then
\begin{equation}\label{eq:strong-loc}
T_{g_k}\xrightarrow{\mathrm{SOT}}cI,
\qquad
H_{g_k}\xrightarrow{\mathrm{SOT}}0,
\qquad
H_{g_k}^*\xrightarrow{\mathrm{SOT}}0,
\end{equation}
and
\begin{equation}\label{eq:comm-loc}
\lim_{k\to\infty}\|[T_{g_k},S_j]\|=0,
\qquad 1\le j\le n.
\end{equation}
Moreover, for every continuous function $h$ on a compact
interval containing $\bigcup_{k\ge1}\sigma(T_{g_k})$,
we have
\begin{equation}\label{eq:calc-loc}
h(T_{g_k})\xrightarrow{\mathrm{SOT}}h(c)I,
\qquad
\|[h(T_{g_k}),S_j]\|\longrightarrow0,
\qquad 1\le j\le n.
\end{equation}
\end{lem}

\begin{proof}
Since $\{g_k\}$ is uniformly bounded in $L^\infty(\partial\mathbb B_n)$
and converges to $c$ in $L^2(\partial\mathbb B_n)$, a standard
density argument shows that $M_{g_k}\to cI$ strongly on
$L^2(\partial\mathbb B_n)$. As the symbols are real-valued,
we also have $M_{g_k}^*\to cI$ strongly.
The assertions in \eqref{eq:strong-loc} now follow by
compressing these multiplication operators to the
corresponding subspaces.

By Lemma~\ref{lem:hankelcompact},
\[
[T_{g_k},S_j]^*=H_{g_k}^*H_{\overline z_j}.
\]
Since $H_{\overline z_j}$ is compact and
$H_{g_k}^*\to 0$ strongly, it follows from
\cite[Chapter~III, Lemma~3.7]{Kato} that
\[
\|[T_{g_k},S_j]\|
=\|H_{g_k}^*H_{\overline z_j}\|
\longrightarrow 0.
\]
Finally, \eqref{eq:calc-loc} follows immediately from
Lemma~\ref{lem:funccalc}.
\end{proof}

\section{Fredholm indices and compressions}
\label{sec:fredholm}
\subsection{Dirac operator and its index}
We recall the exterior algebra construction associated with
a commuting tuple of operators (see
\cite{Arveson,KaadNest}).
Let
\[
\Lambda=\bigoplus_{q=0}^n\bigwedge\nolimits^q\mathbb C^n
\]
be the exterior algebra of $\mathbb C^n$, where
$\bigwedge^0\mathbb C^n=\mathbb C$.
Note that $\dim\Lambda=2^n$.
If $\{e_1,\ldots,e_n\}$ is the standard orthonormal basis
of $\mathbb C^n$, we equip $\Lambda$ with the inner product
for which the vectors
\[
e_{j_1}\wedge\cdots\wedge e_{j_q},
\qquad 1\le j_1<\cdots<j_q\le n,
\quad 0\le q\le n,
\]
form an orthonormal basis, with the empty product
identified with $1$.

The natural $\mathbb Z_2$-grading of $\Lambda$ is given by
$\Lambda=\Lambda^+\oplus\Lambda^-$, where
\[
\Lambda^+=\bigoplus_{q\ \mathrm{even}}
\bigwedge\nolimits^q\mathbb C^n,
\qquad
\Lambda^-=\bigoplus_{q\ \mathrm{odd}}
\bigwedge\nolimits^q\mathbb C^n.
\]
For $1\le j\le n$, let $C_j$ denote the exterior
multiplication operator on $\Lambda$ defined by
\[
C_j\omega=e_j\wedge\omega,
\qquad \omega\in\Lambda.
\]
The operators $C_j$ and their adjoints $C_j^*$ satisfy
the canonical anticommutation relations
\begin{equation}\label{eq:car}
C_iC_j+C_jC_i=0,
\qquad
C_iC_j^*+C_j^*C_i=\delta_{ij}I_\Lambda.
\end{equation}

Set
\[
\mathcal E=H^2(\partial\mathbb B_n)\otimes\Lambda,
\qquad
\mathcal E^\pm=H^2(\partial\mathbb B_n)\otimes\Lambda^\pm.
\]
Let $T=(T_1,\ldots,T_n)$ be a commuting tuple of bounded
operators on $H^2(\partial\mathbb B_n)$.
The associated Koszul differential and Dirac operator
on $\mathcal E$ are defined by
\begin{equation}\label{eq:Ddef}
d_T=\sum_{j=1}^n T_j\otimes C_j,
\qquad D_T=d_T+d_T^*.
\end{equation}
The commutativity of $T_1,\ldots,T_n$, together with
\eqref{eq:car}, implies that $d_T^2=0$.
Moreover, $D_T$ is self-adjoint and odd with respect to
the grading $\mathcal E=\mathcal E^+\oplus\mathcal E^-$.
Hence
\[
D_T=\begin{pmatrix}
0&D_T^-\\
D_T^+&0
\end{pmatrix},
\qquad
D_T^+:\mathcal E^+\longrightarrow\mathcal E^-,
\quad D_T^-=(D_T^+)^*.
\]
In what follows, we take $T=S=(S_1,\ldots,S_n)$
and write $d=d_S$ and $D=D_S$.

\begin{prop}\label{prop:D-index}
The Dirac operator $D$ associated with the spherical shift
$S=(S_1,\ldots,S_n)$ satisfies
\begin{equation}\label{eq:Desunitary}
D^2-I_{\mathcal E}\in\mathcal K(\mathcal E).
\end{equation}
In particular, $D^+$ is Fredholm and
\begin{equation}\label{eq:indexone}
    \operatorname{ind}D^+=(-1)^n.
\end{equation}
\end{prop}

\begin{proof}
Since $d^2=0$, the anticommutation relations \eqref{eq:car}
yield
\begin{equation*}
\begin{aligned}
D^2
&=dd^*+d^*d\\
&=\sum_{i=1}^n S_iS_i^*\otimes I_\Lambda
+\sum_{i,j=1}^n
[S_j^*,S_i]\otimes C_j^*C_i.
\end{aligned}
\end{equation*}
By \eqref{eq:essential-shift} and the finite dimensionality
of $\Lambda$, we obtain \eqref{eq:Desunitary}.
Consequently,
\[
D^-D^+-I_{\mathcal E^+}\in\mathcal K(\mathcal E^+),
\qquad
D^+D^--I_{\mathcal E^-}\in\mathcal K(\mathcal E^-).
\]
Thus $D^+$ is Fredholm by Atkinson's theorem
(see \cite[Section~XI.2]{Conway}).

To compute its index, for each $N\ge0$, set
\[
\mathcal E_N
=\bigoplus_{q=0}^n
H_{N-n+q}\otimes\bigwedge\nolimits^q\mathbb C^n,
\]
where $H_m=\{0\}$ for $m<0$. Then
\[
\mathcal E=\bigoplus_{N=0}^\infty\mathcal E_N.
\]
By construction, each $\mathcal E_N$ is invariant under both
$d$ and $d^*$, and hence reduces $D$.

Put
\[
\mathcal E_N^\pm=\mathcal E_N\cap\mathcal E^\pm,
\qquad
D_N^+=D^+|_{\mathcal E_N^+}:
\mathcal E_N^+\longrightarrow\mathcal E_N^-.
\]
We have
\[
D^+=\bigoplus_{N=0}^\infty D_N^+.
\]
Since $D^+$ is Fredholm, both $\ker D^+$ and $\ker (D^+)^*$
are finite dimensional. The corresponding orthogonal
decompositions of these kernels yield
\[
\operatorname{ind}D^+
=\sum_{N=0}^\infty\operatorname{ind}D_N^+,
\]
where only finitely many summands are nonzero.

Since $\mathcal E_N^\pm$ are finite dimensional,
\[
\operatorname{ind}D_N^+
=\dim\mathcal E_N^+-\dim\mathcal E_N^-
=\sum_{q=0}^n(-1)^q\binom nq d_{N-n+q},
\]
where
\[
d_m=\dim H_m=\binom{m+n-1}{n-1}
\quad(m\ge0),
\qquad d_m=0\quad(m<0).
\]
For $|t|<1$, we therefore have
\[
\begin{aligned}
\sum_{N=0}^\infty\operatorname{ind}D_N^+\,t^N
&=\left(\sum_{q=0}^n(-1)^q\binom nq t^{n-q}\right)
\left(\sum_{m=0}^\infty d_m t^m\right)\\
&=(t-1)^n(1-t)^{-n}\\
&=(-1)^n.
\end{aligned}
\]
Comparing coefficients gives
$\operatorname{ind}D_0^+=(-1)^n$ and
$\operatorname{ind}D_N^+=0$ for every $N\ge1$.
This proves \eqref{eq:indexone}.
\end{proof}

\subsection{Fredholm indices of compressions}
Let $Q$ be an orthogonal projection on
$H^2(\partial\mathbb B_n)$ and set
$\widehat Q=Q\otimes I_\Lambda$.
Since $\widehat Q$ preserves the grading of $\mathcal E$,
the compression
\[
D_Q=\widehat QD\big|_{\widehat Q\mathcal E}
\]
is an odd self-adjoint operator on $\widehat Q\mathcal E$.
Denote its component by
\[
D_Q^+:\widehat Q\mathcal E^+
\longrightarrow\widehat Q\mathcal E^-.
\]
We establish the following additivity formula for the
Fredholm index.

\begin{prop}\label{prop:additivity}
Let $Q$ be an orthogonal projection on
$H^2(\partial\mathbb B_n)$ such that $[Q,S_j]$ is compact
for every $1\le j\le n$.
Then both $D_Q^+$ and $D_{I-Q}^+$ are Fredholm, and
\begin{equation}\label{eq:additivity}
\operatorname{ind}D_Q^+
+\operatorname{ind}D_{I-Q}^+
=\operatorname{ind}D^+=(-1)^n.
\end{equation}
\end{prop}

\begin{proof}
Since $Q=Q^*$, the compactness of $[Q,S_j]$ implies
that of $[Q,S_j^*]$. Hence $[\widehat Q,D]$ is compact.

On $\widehat Q\mathcal E$, we have
\begin{equation}\label{eq:squarecompression}
\begin{aligned}
D_Q^2-I_{\widehat Q\mathcal E}
=\widehat Q(D^2-I_{\mathcal E})
   \big|_{\widehat Q\mathcal E}
   -\widehat QD(I-\widehat Q)D
   \big|_{\widehat Q\mathcal E}.
\end{aligned}
\end{equation}
The first term on the right hand side is compact by
Proposition~\ref{prop:D-index}, while the second is compact
since
\[
\widehat QD(I-\widehat Q)
=[\widehat Q,D](I-\widehat Q).
\]
Thus $D_Q^2-I_{\widehat Q\mathcal E}$ is compact.
Consequently,
\[
D_Q^-D_Q^+-I_{\widehat Q\mathcal E^+}
\in\mathcal K(\widehat Q\mathcal E^+)
\]
and
\[
D_Q^+D_Q^--I_{\widehat Q\mathcal E^-}
\in\mathcal K(\widehat Q\mathcal E^-).
\]
By Atkinson's theorem, $D_Q^+$ is Fredholm.
Replacing $Q$ by $I-Q$ gives the Fredholmness of
$D_{I-Q}^+$.

With respect to the orthogonal decompositions
\[
\mathcal E^\pm
=\widehat Q\mathcal E^\pm
\oplus(I-\widehat Q)\mathcal E^\pm,
\]
the diagonal blocks of $D^+$ are precisely $D_Q^+$
and $D_{I-Q}^+$, while the off-diagonal blocks are
restrictions of $\widehat QD(I-\widehat Q)$ and
$(I-\widehat Q)D\widehat Q$, respectively.
Since $[\widehat Q,D]$ is compact, both off-diagonal
blocks are compact. Consequently,
\[
D^+-(D_Q^+\oplus D_{I-Q}^+)
\in\mathcal K(\mathcal E^+,\mathcal E^-).
\]
The additivity of the Fredholm index under direct sums
and its invariance under compact perturbations yield
\eqref{eq:additivity}.
\end{proof}

The preceding index formula
yields the following limiting result.

\begin{lem}\label{lem:indexlimits}
Let $\{Q_k\}$ be a sequence of orthogonal projections on
$H^2(\partial\mathbb B_n)$ such that $[Q_k,S_j]$ is compact
for every $k$ and $1\le j\le n$. Suppose that
\[
\|[\widehat Q_k,D]\|\longrightarrow0.
\]
Then the following assertions hold.
\begin{enumerate}
\item[(i)] If $Q_k\to 0$ strongly, then $D_{Q_k}$
is invertible and $\operatorname{ind}D_{Q_k}^+=0$
for all sufficiently large $k$.

\item[(ii)] If $Q_k\to I$ strongly, then
$\operatorname{ind}D_{Q_k}^+=(-1)^n$
for all sufficiently large $k$.
\end{enumerate}
\end{lem}

\begin{proof}
(i) Since $\Lambda$ is finite dimensional, the strong
convergence $Q_k\to 0$ implies that
$\widehat Q_k\to 0$ strongly on $\mathcal E$.
By Proposition~\ref{prop:D-index},
$D^2-I_{\mathcal E}$ is compact. Consequently,
\[
\|\widehat Q_k(D^2-I_{\mathcal E})\widehat Q_k\|
\longrightarrow 0.
\]
Since $D$ is self-adjoint, it follows from
\eqref{eq:squarecompression} that
\begin{equation*}
\begin{aligned}
\|D_{Q_k}^2-I_{\widehat Q_k\mathcal E}\|
&\le
\|\widehat Q_k(D^2-I_{\mathcal E})\widehat Q_k\|
+\|\widehat Q_kD(I-\widehat Q_k)\|^2\\
&\le
\|\widehat Q_k(D^2-I_{\mathcal E})\widehat Q_k\|
+\|[\widehat Q_k,D]\|^2\\
&\longrightarrow0.
\end{aligned}
\end{equation*}
It follows that $D_{Q_k}^2$, and hence $D_{Q_k}$,
is invertible for all sufficiently large $k$.
Since $D_{Q_k}$ is odd and self-adjoint,
the restriction $D_{Q_k}^+$ is also invertible. In particular,
\[
\operatorname{ind}D_{Q_k}^+=0
\]
for all sufficiently large $k$.

\smallskip
\noindent
(ii) Suppose that $Q_k\to I$ strongly.
Then $I-Q_k\to0$ strongly and
\[
[I-\widehat Q_k,D]=-[\widehat Q_k,D].
\]
Applying part~(i) to the complementary projections
$I-Q_k$, we obtain
$\operatorname{ind}D_{I-Q_k}^+=0$
for all sufficiently large $k$.
The index formula \eqref{eq:additivity} now yields
\[
\operatorname{ind}D_{Q_k}^+
=\operatorname{ind}D^+=(-1)^n
\]
for all sufficiently large $k$.
\end{proof}

We next establish the invariance of the Fredholm index
under the unitary transformations induced by
automorphisms of the unit ball.

\begin{prop}\label{prop:invariance}
Let $g\in L^\infty(\partial\mathbb B_n)$ be real-valued,
and let $h$ be a continuous function on an
interval containing $\sigma(T_g)$ such that $Q=h(T_g)$
is an orthogonal projection. For $a\in\mathbb B_n$, set
\begin{equation*}
Q_a=U_aQU_a^*=h(T_{g\circ\varphi_a}).
\end{equation*}
Then $D_{Q_a}^+$ is Fredholm and
\[
\operatorname{ind}D_{Q_a}^+
=\operatorname{ind}D_Q^+
\qquad(a\in\mathbb B_n).
\]
\end{prop}

\begin{proof}
By Lemmas~\ref{lem:hankelcompact} and \ref{lem:funccalc},
$[Q_a,S_j]$ is compact for every $a\in\mathbb B_n$
and $1\le j\le n$.
Proposition~\ref{prop:additivity} therefore implies
that $D_{Q_a}^+$ is Fredholm.

Set $\widehat U_a=U_a\otimes I_\Lambda$.
Since $\widehat U_a$ preserves the grading of
$\mathcal E$ and satisfies
\[
\widehat Q_a=\widehat U_a\widehat Q\widehat U_a^*,
\]
its restrictions to $\widehat Q\mathcal E^\pm$ are
unitary operators,
respectively.
We may thus regard the compressions $D_{Q_a}$,
up to unitary equivalence, as operators on the fixed
Hilbert space $\widehat Q\mathcal E$. More precisely,
\begin{equation}\label{eq:fixedconj}
\widehat U_a^*D_{Q_a}\widehat U_a
=\widehat QD(a)\big|_{\widehat Q\mathcal E},
\end{equation}
where $D(a)=\widehat U_a^*D\widehat U_a$.
Since $U_a=U_a^*$, the covariance formula
\eqref{eq:covariance} gives
\begin{equation}\label{eq:Da}
\begin{aligned}
D(a)=\sum_{j=1}^n\bigl(
&M_{\varphi_{a,j}}|_{H^2(\partial\mathbb B_n)}
\otimes C_j\\
&+(M_{\varphi_{a,j}}|_{H^2(\partial\mathbb B_n)})^*
\otimes C_j^*\bigr).
\end{aligned}
\end{equation}
By Lemma~\ref{lem:covariance}, the mapping
$a\mapsto D(a)$ is continuous in the operator norm.
Consequently, the family of operators
\[
F_a=\widehat QD(a)\big|_{\widehat Q\mathcal E^+}:
\widehat Q\mathcal E^+\longrightarrow
\widehat Q\mathcal E^-,
\qquad a\in\mathbb B_n,
\]
is norm-continuous. By \eqref{eq:fixedconj}, each $F_a$
is Fredholm and
$
\operatorname{ind}F_a=\operatorname{ind}D_{Q_a}^+.
$
By the continuity of the Fredholm index
\cite[Theorem~5.36]{Douglas} and the connectedness
of $\mathbb B_n$, we obtain
\[
\operatorname{ind}F_a=\operatorname{ind}F_0
\qquad(a\in\mathbb B_n).
\]
Noting that $\varphi_0(z)=-z$, we obtain $D(0)=-D$
from \eqref{eq:Da}, and hence $F_0=-D_Q^+$.
It follows that
\[
\operatorname{ind}D_{Q_a}^+
=\operatorname{ind}F_a
=\operatorname{ind}F_0
=\operatorname{ind}D_Q^+.
\]
This completes the proof.
\end{proof}

\section{Proof of the main theorem}
\label{sec:mainproof}

\begin{proof}[Proof of Theorem~\ref{thm}]
Set
\[
m=\operatorname*{ess\,inf}_{\partial\mathbb B_n} f,
\qquad
M=\operatorname*{ess\,sup}_{\partial\mathbb B_n} f.
\]
The operator inequalities $mI\le T_f\le MI$ imply that
$\sigma(T_f)\subseteq[m,M]$.
If $m=M$, then $T_f=mI$, and the assertion is immediate.
We therefore assume that $m<M$.

Suppose, towards a contradiction, that there exists
$\lambda\in(m,M)\setminus\sigma(T_f)$.
Choose
\[
0<\varepsilon<\frac12
\min\{\lambda-m,M-\lambda,
\operatorname{dist}(\lambda,\sigma(T_f))\}.
\]
In particular,
\[
\sigma(T_f)\cap
[\lambda-2\varepsilon,\lambda+2\varepsilon]
=\varnothing,
\]
and the sets
\[
E_-=\{\zeta\in\partial\mathbb B_n:
f(\zeta)<\lambda-2\varepsilon\},
\qquad
E_+=\{\zeta\in\partial\mathbb B_n:
f(\zeta)>\lambda+2\varepsilon\}
\]
both have positive surface measure.

Let $h:\mathbb R\to[0,1]$ be a continuous function such that
\[
h(t)=0\quad\text{for }t\le\lambda-\varepsilon,
\qquad
h(t)=1\quad\text{for }t\ge\lambda+\varepsilon.
\]
Since $h^2=h$ on $\sigma(T_f)$, the continuous functional
calculus shows that $Q=h(T_f)$ is an orthogonal projection.
Moreover, by Lemmas~\ref{lem:hankelcompact}
and \ref{lem:funccalc},
\[
[Q,S_j]\in\mathcal K(H^2(\partial\mathbb B_n)),
\qquad 1\le j\le n.
\]
For each $a\in\mathbb B_n$, set $Q_a=U_aQU_a^*$.
By Proposition~\ref{prop:invariance},
\begin{equation}\label{eq:nuzero}
\operatorname{ind}D_{Q_a}^+
=\operatorname{ind}D_Q^+,
\qquad a\in\mathbb B_n.
\end{equation}

Since $E_-$ and $E_+$ have positive surface measure,
by Lemma~\ref{lem:scalar} we may choose
$\zeta_\pm\in E_\pm$ such that
\[
f\circ\varphi_{r\zeta_\pm}
\longrightarrow f(\zeta_\pm)
\quad\text{in }L^2(\partial\mathbb B_n)
\quad\text{as }r\to1^-.
\]
Fix a sequence $\{r_k\}\subset(0,1)$ converging to $1$
and set $a_k^\pm=r_k\zeta_\pm$.
By Lemma~\ref{lem:covariance} and the continuous
functional calculus,
\[
Q_{a_k^\pm}
=h\big(T_{f\circ\varphi_{a_k^\pm}}\big).
\]
Since the transformed symbols are uniformly bounded
in $L^\infty(\partial\mathbb B_n)$,
Lemma~\ref{lem:localcomm} gives
\[
Q_{a_k^\pm}\xrightarrow{\mathrm{SOT}}
h(f(\zeta_\pm))I,
\qquad
\|[Q_{a_k^\pm},S_j]\|\longrightarrow0,
\quad 1\le j\le n.
\]
Moreover, by the definition of $D$,
\[
\|[\widehat Q_{a_k^\pm},D]\|\longrightarrow0.
\]
By the choice of $E_\pm$ and $h$, we have
\[
Q_{a_k^-}\xrightarrow{\mathrm{SOT}}0,
\qquad
Q_{a_k^+}\xrightarrow{\mathrm{SOT}}I.
\]
It follows from Lemma~\ref{lem:indexlimits} that,
for all sufficiently large $k$,
\[
\operatorname{ind}D_{Q_{a_k^-}}^+=0,
\qquad
\operatorname{ind}D_{Q_{a_k^+}}^+
=\operatorname{ind}D^+.
\]
Combining these identities with \eqref{eq:nuzero},
we obtain
\[
0=\operatorname{ind}D_Q^+
=\operatorname{ind}D^+=(-1)^n,
\]
a contradiction.

Thus $(m,M)\subseteq\sigma(T_f)$.
Since $\sigma(T_f)$ is closed and contained in $[m,M]$,
we conclude that
\[
\sigma(T_f)=[m,M].
\]
Finally, since $T_f$ is self-adjoint and
$\sigma(T_f)=[m,M]$ has no isolated points,
it follows from \cite[Chapter~XI,
Proposition~4.6]{Conway} that
\[
\sigma_{\mathrm e}(T_f)
=\sigma(T_f)=[m,M].
\]
\end{proof}

\section{Further consequences}
\label{sec:consequences}

We first apply  Theorem~\ref{thm} to Toeplitz operators with
characteristic-function symbols and obtain an exact
description of the spectral properties of the associated
Hankel operators. 

\begin{cor}\label{cor:char}
Let $E\subseteq\partial\mathbb B_n$ be a measurable set
with $0<\sigma(E)<1$, and let $\mathbf1_E$ denote its
characteristic function. Then
\begin{equation}\label{eq:char}
\sigma(T_{\mathbf1_E})
=\sigma_{\mathrm e}(T_{\mathbf1_E})=[0,1].
\end{equation}
Moreover,
\[
\sigma(H_{\mathbf1_E}^*H_{\mathbf1_E})
=\sigma_{\mathrm e}(H_{\mathbf1_E}^*H_{\mathbf1_E})
=\left[0,\frac14\right].
\]
In particular,
\[
\|H_{\mathbf1_E}\|
=\|H_{\mathbf1_E}\|_{\mathrm e}
=\frac12.
\]
\end{cor}

\begin{proof}
The spectral identity \eqref{eq:char} follows immediately
from  Theorem~\ref{thm}. Since $\mathbf1_E$ is real-valued and
idempotent, the Toeplitz--Hankel identity gives
\[
H_{\mathbf1_E}^*H_{\mathbf1_E}
=T_{\mathbf1_E}-T_{\mathbf1_E}^2.
\]
By \eqref{eq:char} and the spectral mapping theorem
for both the spectrum and the essential spectrum, we have
\[
\sigma(H_{\mathbf1_E}^*H_{\mathbf1_E})
=\sigma_{\mathrm e}(H_{\mathbf1_E}^*H_{\mathbf1_E})
=\left[0,\frac14\right].
\]
Consequently,
\[
\|H_{\mathbf1_E}\|^2
=\|H_{\mathbf1_E}^*H_{\mathbf1_E}\|
=\|H_{\mathbf1_E}^*H_{\mathbf1_E}\|_{\mathrm e}
=\|H_{\mathbf1_E}\|_{\mathrm e}^2
=\frac14.
\]
This completes the proof.
\end{proof}

We next show that the spectral conclusion of  Theorem~\ref{thm}
remains valid on weighted Hardy spaces whose norms are
equivalent to the usual Hardy norm.
Given a weight $w$ satisfying
$0<c\le w\le C<\infty$ almost everywhere on
$\partial\mathbb B_n$, let $H^2(w\,d\sigma)$ denote
the closure of the holomorphic polynomials in
$L^2(\partial\mathbb B_n,w\,d\sigma)$, and let $P_w$
be the corresponding orthogonal projection.
For $f\in L^\infty(\partial\mathbb B_n)$, define
\[
T_f^{(w)}
=P_wM_f\big|_{H^2(w\,d\sigma)},
\qquad
H_f^{(w)}
=(I-P_w)M_f\big|_{H^2(w\,d\sigma)}.
\]

\begin{cor}\label{cor:weighted}
Let $w\in L^\infty(\partial\mathbb B_n)$ satisfy
$0<c\le w\le C<\infty$ almost everywhere.
Then, for every real-valued $f\in L^\infty(\partial\mathbb B_n)$,
\begin{equation}\label{eq:weightedB}
\sigma(T_f^{(w)})
=\sigma_{\mathrm e}(T_f^{(w)})
=\left[
\operatorname*{ess\,inf}_{\partial\mathbb B_n}f,
\operatorname*{ess\,sup}_{\partial\mathbb B_n}f
\right].
\end{equation}
In particular, for every measurable set
$E\subseteq\partial\mathbb B_n$ with $0<\sigma(E)<1$,
we have
\[
\sigma\big((H_{\mathbf1_E}^{(w)})^*
H_{\mathbf1_E}^{(w)}\big)
=\sigma_{\mathrm e}\big((H_{\mathbf1_E}^{(w)})^*
H_{\mathbf1_E}^{(w)}\big)
=\left[0,\frac14\right]
\]
and
\[
\|H_{\mathbf1_E}^{(w)}\|
=\|H_{\mathbf1_E}^{(w)}\|_{\mathrm e}
=\frac12.
\]
\end{cor}

\begin{proof}
Since $0<c\le w\le C<\infty$ almost everywhere,
the weighted and unweighted Hardy norms are equivalent.
Moreover, the inequalities $cI\le T_w\le CI$
imply that $T_w$ is invertible.

Let $\langle\cdot,\cdot\rangle_w$ denote the weighted
inner product. For $h,k\in H^2(\partial\mathbb B_n)$,
we have
\[
\langle T_wT_f^{(w)}h,k\rangle
=\langle T_f^{(w)}h,k\rangle_w
=\langle fh,k\rangle_w
=\langle T_{wf}h,k\rangle.
\]
Consequently,
\begin{equation}\label{eq:weightedidentity}
T_w(T_f^{(w)}-\lambda I)=T_{w(f-\lambda)},
\qquad \lambda\in\mathbb C.
\end{equation}
By the invertibility of $T_w$ and the equivalence of norms,
$T_f^{(w)}-\lambda I$ is Fredholm if and only if
$T_{w(f-\lambda)}$ is Fredholm.

Set
\[
m=\operatorname*{ess\,inf}_{\partial\mathbb B_n}f,
\qquad
M=\operatorname*{ess\,sup}_{\partial\mathbb B_n}f.
\]
The operator inequalities $mI\le T_f^{(w)}\le MI$
on $H^2(w\,d\sigma)$ yield
$\sigma(T_f^{(w)})\subseteq[m,M]$.
If $m=M$, then $T_f^{(w)}=mI$, and the assertion follows.
Suppose that $m<M$ and let $\lambda\in(m,M)$.
Since $w\ge c>0$ almost everywhere, we have
\[
\operatorname*{ess\,inf}w(f-\lambda)
<0<
\operatorname*{ess\,sup}w(f-\lambda).
\]
By Theorem~\ref{thm} and \eqref{eq:weightedidentity},
we obtain
\[
0\in\sigma_{\mathrm e}(T_{w(f-\lambda)}),
\qquad
\lambda\in\sigma_{\mathrm e}(T_f^{(w)}).
\]
Hence, by the closedness of the essential spectrum,
\[
[m,M]\subseteq\sigma_{\mathrm e}(T_f^{(w)})
\subseteq\sigma(T_f^{(w)})\subseteq[m,M],
\]
which proves \eqref{eq:weightedB}.

The remaining assertions follow from \eqref{eq:weightedB}
by the same argument as in the proof of
Corollary~\ref{cor:char}.
\end{proof}

\medskip
\noindent\textbf{AI Declaration.}
During the preparation of this manuscript, the authors
used ChatGPT (OpenAI) to assist with language refinement,
\LaTeX{} formatting, and checking mathematical arguments for potential gaps.
All mathematical statements, proofs, and conclusions
were independently verified by the authors, who
take full responsibility for the content of this article.

\medskip
\noindent\textbf{Acknowledgments.}
The authors are grateful to Professor Jingbo Xia for
his careful reading of the manuscript and valuable
suggestions for its improvement.
T. Yu was supported by the National Natural Science
Foundation of China (Grant No.~12571132).

\end{document}